\documentclass[11pt]{article}
\usepackage{epic,latexsym}
\usepackage{color}
\usepackage{tikz}
\usepackage{comment}
\usepackage{amssymb}
\usepackage{mathtools}
\usepackage{amsthm}
\usepackage{url}

\newtheorem{theorem}{Theorem}[section]

\newtheorem{proposition}[theorem]{Proposition}

\newtheorem{corollary}[theorem]{Corollary}

\newtheorem{lemma}[theorem]{Lemma}

\newcommand{\Dim}{{\rm Dim}}
\newcommand{\ex}{{\rm ex}}
\newcommand{\QED}{$\Box$}

\newcommand{\dimS}{\dim_{\rm S}}
\newcommand{\DimS}{\Dim_{\rm S}}

\newcommand{\modo}{{\rm mod \,}}
\newcommand{\2}{ \vspace{0.2cm} }
\newcommand{\1}{ \vspace{0.1cm} }

\let\oldenumerate\enumerate
\renewcommand{\enumerate}{
  \oldenumerate
  \setlength{\itemsep}{0pt}
  \setlength{\parskip}{0pt}
  \setlength{\parsep}{0pt}
}

\begin{document}

\title{Resolvability and convexity properties in the Sierpi\'{n}ski product of graphs}

\author{$^1$Michael A. Henning \qquad
$^{2,3,4}$Sandi Klav\v zar\qquad $^5$Ismael G. Yero
\\ \\
\small $^1$Department of Mathematics and Applied Mathematics \\
\small University of Johannesburg \\
\small Auckland Park, 2006 South Africa\\
\small \tt mahenning@uj.ac.za
\\ \\
\small $^2$Faculty of Mathematics and Physics \\
\small University of Ljubljana, Slovenia\\
\small \tt sandi.klavzar@fmf.uni-lj.si\\
\\
\small $^3$Institute of Mathematics, Physics and Mechanics \\
\small Ljubljana, Slovenia\\
\\
\small $^4$ Faculty of Natural Sciences and Mathematics \\
\small University of Maribor, Slovenia \\
\\
\small $^5$Departamento de Matem\'aticas\\
\small Universidad de C\'adiz, Algeciras, Spain \\
\small \tt ismael.gonzalez@uca.es\\
}
\date{}
\maketitle

\begin{abstract}
Let $G$ and $H$ be graphs and let $f \colon V(G)\rightarrow V(H)$ be a function. The Sierpi\'{n}ski product of $G$ and $H$ with respect to $f$, denoted by  $G \otimes _f H$, is defined as the graph on the vertex set $V(G)\times V(H)$, consisting of $|V(G)|$ copies of $H$; for every edge $gg'$ of $G$ there is an edge between copies $gH$ and $g'H$ of $H$ associated with the vertices $g$ and $g'$ of $G$, respectively, of the form $(g,f(g'))(g',f(g))$. The Sierpi\'{n}ski metric dimension and the upper Sierpi\'{n}ski metric dimension of two graphs are determined. Closed formulas are determined for Sierpi\'{n}ski products of trees, and for Sierpi\'{n}ski products of two cycles where the second factor is a triangle. We also prove that the layers with respect to the second factor in a Sierpi\'{n}ski product graph are convex.
\end{abstract}


\newpage

\section{Introduction}

Sierpi\'{n}ski graphs represent a very interesting and widely studied family of graphs. They were introduced in 1997 in the paper~\cite{klavzar-1997}, where the primary motivation for their introduction was the intrinsic link to the Tower of Hanoi problem, for the latter problem see the book~\cite{hinz-2018}. Intensive research of Sierpi\'{n}ski graphs led to a review article~\cite{hinz-2017} in which state of the art up to 2017 is summarized and unified approach to Sierpi\'{n}ski-type graph families is also proposed. Later research on Sierpi\'{n}ski graphs includes~\cite{bresar-2018, estrada-moreno-2018, estrada-moreno-2019, farrokhi-2021, khatibi-2020, liu-2021, varghese-2021}.

In this paper we study a recent generalization of Sierpi\'{n}ski graphs proposed by Kovi\v{c}, Pisanski, Zemlji\v{c}, and \v{Z}itnik in~\cite{kpzz-2022}. Let $G$ and $H$ be graphs and let $f \colon V(G)\rightarrow V(H)$ be an arbitrary function. The \textit{Sierpi\'{n}ski product of graphs $G$ and $H$ with respect to $f$}, denoted by  $G \otimes _f H$, is defined as the graph on the vertex set $V(G)\times V(H)$ with edges of two types:
\begin{itemize}
    \item \emph{Type-$1$ edge}:
    $(g,h)(g,h')$ is an edge of  $G \otimes _f H$ for every vertex $g\in V(G)$ and every edge $hh' \in E(H)$,
    \item \emph{Type-$2$ edge}: $(g,f(g'))(g',f(g))$ is an edge of $G \otimes _f H$ for every edge $gg' \in E(G)$.
\end{itemize}

We observe that the edges of Type-$1$ induce $n(G) = |V(G)|$ copies of the graph $H$ in the Sierpi\'{n}ski product $G \otimes _f H$. For each vertex $g \in V(G)$, we let $gH$ be the copy of $H$ corresponding to the vertex $g$. A Type-$2$ edge joins vertices from different copies of $H$ in $G \otimes _f H$, and is called a \emph{connecting edge} of $G \otimes _f H$. A vertex incident with a connecting edge is called a \emph{connecting vertex}. We observe that two different copies of $H$ in $G \otimes _f H$ are joined by at most one edge. We denote by $H^G$ be the family of functions from $V(G)$ to $V(H)$.

It might be readily observed that the Sierpi\'{n}ski product is closely related to other product graphs. For instance, by considering a constant function $f$ in the product, we obtain graphs which are indeed the same as the so-called rooted product graphs (see \cite{godsil-1978} for its definition). Also, selecting the identity function ${\rm id}\in G^G$, the Sierpi\'{n}ski product $G \otimes_{\rm id} G$ is the (first iteration of the) generalized Sierpi\'{n}ski graph in the sense of~\cite{gravier-2011}. Moreover, a Sierpi\'{n}ski product can also be considered as a subgraph of the (Cartesian, strong or lexicographic) product. Consequently, any contribution to the study of the Sierpi\'{n}ski product could give some more knowledge on these related products.

In the next two subsections we give motivation,  basic terminology, and notation concerning the classical metric dimension of graphs, and introduce the study of the Sierpi\'{n}ski metric dimension and the upper Sierpi\'{n}ski metric dimension. Thereafter in Section~\ref{S:trees} we determine the upper Sierpi\'{n}ski metric dimension for Sierpi\'{n}ski products of arbitrary trees. A general lower bound is established for the Sierpi\'{n}ski metric dimension for products of two trees, and an exact formula when the first factor is a path. In Section~\ref{S:cycles} a closed formula is determined for both dimensions when the first factor in the product is an arbitrary cycle and the second factor a triangle. In Section~\ref{S:convex} we prove that the layers with respect to the second factor in a Sierpi\'{n}ski product graph are convex. In Section~\ref{S:remark} we pose several open problems.

\subsection{The metric dimension of graphs}
\label{S:metric-dim}

The \emph{distance} between two vertices $u$ and $v$ in a connected graph $G$, denoted $d_G(u,v)$, is the number of edges in a shortest path from $u$ to $v$, that is, $d_G(u,v)$ is the minimum length of a $u,v$-path in $G$. Given an ordered subset $S = \{v_1,\ldots,v_k\}$ of vertices in $G$, the \emph{metric} $S$-\emph{representation} of a vertex $v$ in $G$ is the $k$-tuple vector
\[
r_G(v|S)=(d_G (v,v_1) ,\cdots, d_G (v,v_k )).
 \]
If every two distinct vertices of $G$ have different metric $S$-representations, then the set $S$ is called a \emph{resolving set} of $G$ (also called a \emph{metric generator}). The \emph{metric dimension} of $G$, denoted by $\dim(G)$, is the cardinality of a smallest possible resolving set in~$G$. A \emph{metric basis} of $G$ is a resolving set of cardinality $\dim(G)$. A vertex $v$ in a graph $G$ is said to \emph{distinguish} (or \emph{resolve}) two vertices $x$ and $y$ if $d_G(v,x) \ne d_G(v,y)$.

The concept of the metric dimension of a graph was birthed independently by Harary and Melter~\cite{HaMe-76} in 1976 and by Slater~\cite{Slater-75} in 1975, and is now well studied in graph theory. To date MathSciNet lists over 380 papers on metric dimension in graphs, covering a large number of different investigations dealing with theoretical and applied results on such parameter.

According to the structural properties of resolving sets in graphs, they can easily be used to model several practical situations in which uniquely recognizing points or locations is required. That was precisely one of the motivations of the seminal works~\cite{HaMe-76} and \cite{Slater-75}, where resolving sets appeared to be used for the location of intruders in networks. Further on, some other related models and applications have appeared here and there. Among them, we remark the recent work \cite{till-2019}, where the authors presented a connection between some metric dimension parameter and the representation of genomic sequences. Among the theoretical studies on this topic, the literature contains a wide range of different contributions, some recent and remarkable articles are for instance \cite{claverol-2021,geneson-2022,Sedlar-2021}. For more information on investigations on the classical version we suggest the fairly complete survey \cite{till-2022+}.

With respect to the theoretical studies, the metric dimension of graph products and graph operations has attracted the attention of several investigations. In this sense, we mention a few interesting contributions related with this exposition due to the relationship between the Sierpi\'{n}ski product and some other products previously mentioned. The metric dimension of Cartesian product graphs has been considered in several works like \cite{Caceres-2007}, for the general case, and among other ones, in \cite{Chau,Jiang2019,Khuller} for some particular examples of Cartesian products. The lexicographic product of graphs has been studied with respect to its metric dimension in \cite{Jannesari,Saputro}, while the strong product has been considered in \cite{Adar,rodriguez-2015}. On the other hand, the metric dimension of the rooted product has been dealt with in \cite{feng-2013,kuziak-2017}.

\subsection{Sierpi\'{n}ski metric dimension}
\label{S:defn}

Let $G$ and $H$ be graphs and $H^G$ be the family of functions from $V(G)$ to $V(H)$. We introduce new types of metric dimension, the \textit{Sierpi\'{n}ski metric dimension}, denoted by $\dimS(G,H)$, as the minimum over all functions $f$ from $H^G$ of the metric dimension of the Sierpi\'{n}ski product with respect to $f$, and \textit{upper Sierpi\'{n}ski metric dimension}, denoted by $\DimS(G, H)$, as the maximum over all functions $f\in H^G$ of the metric dimension of the Sierpi\'{n}ski product with respect to $f$. That is,
\[
\dimS(G, H) \coloneqq \min_{f\in H^G}\{\dim(G\otimes _f H)\}
\]
and
\[
\DimS(G, H) \coloneqq \max _{f\in H^G}\{\dim(G\otimes _f H)\}\,.
\]

We might remark that the classical metric dimension of Sierpi\'{n}ski graphs was already studied in \cite{klavzar-2018}, as well as, that of the generalized Sierpi\'nski graphs over stars was considered in \cite{alizadeh-2019}.

\section{Sierpi\'nski products of trees}
\label{S:trees}

A vertex of degree at least~$3$ in a tree $T$ is called a \emph{branch vertex} (also called a \emph{major vertex} in the literature). A leaf $u$ of $T$ is called a \emph{terminal leaf} of a branch vertex $v$ of $T$ if $d_T(u,v) < d_T(u,w)$ for every other branch vertex $w$ of $T$. The \emph{terminal degree} of a branch vertex $v$ is the number of terminal leaves associated with~$v$. A branch vertex $v$ of $T$ is an \emph{exterior branch vertex} of $T$ if it has positive terminal degree. The path from a terminal leaf to the vertex immediately preceding the branch vertex that it is closest to is called a \emph{terminal path}. Thus, every vertex on a terminal path in $T$ is either a leaf of $T$ or has degree~$2$ in $T$. A vertex on a terminal path that has degree~$2$ in $T$ is called an \emph{internal terminal vertex}. Equivalently, every vertex on a terminal path that is not a terminal leaf, is an internal terminal vertex. Thus if $u$ is an internal terminal vertex in $T$, then the vertex $u$ is an internal vertex of a path $P$ that joins a leaf and a branch vertex closest to that leaf in $T$ where every internal vertex of $P$ has degree~$2$ in $T$.

Let $n_1(T)$ denote the number of leaves of $T$, and let $\ex(T)$ denote the number of exterior branch vertices of $T$. The formula for the metric dimension of a tree reads as follows.

\begin{theorem}{\rm (\cite{HaMe-76,Slater-75})}
\label{thm:mdim-tree}
If $T$ is a tree that is not a path, then
\begin{equation}
\label{bounDimUsandoExterior}
\dim(T) = n_1(T) - \ex(T).
\end{equation}
\end{theorem}

It is clear that $\dim(P_n) = 1$. Combining this fact with Theorem~\ref{thm:mdim-tree} yields the following consequence.

\begin{corollary}
\label{known-mdim-tree}
If $T$ is a tree, then $\dim(T) = 1$ if $T$ is a path and $\dim(T) \ge 2$ if $T$ is not a path.
\end{corollary}

Let $T$ be a tree that is not a path, and let $v_1, \ldots, v_k$ be the exterior branch vertices in $T$ that have terminal degree at least~$2$. If the exterior branch vertex $v_i$ has terminal degree~$\ell_i \ge 2$ and if $L_i$ is a set consisting of all terminal leaves but one associated with $v_i$ for all $i \in [k]$, then~\eqref{bounDimUsandoExterior} can be equivalently stated as:
\[
\dim(T) = \sum_{i=1}^k (\ell_i - 1),
\]
and the set
\[
B(T) = \bigcup_{i=1}^k L_i \2
\]
is a metric basis of $T$ (of cardinality~$\dim(T)$). We call the basis $B(T)$ a \emph{standard metric basis} of $T$. Thus, every vertex in a standard metric basis of a tree $T$ is a leaf, and such a basis contains all but one selected fixed leaf associated with the exterior branch vertex of terminal degree at least~$2$ in $T$.

\subsection{Upper Sierpi\'{n}ski metric dimension in trees}

In this section we determine the upper Sierpi\'{n}ski metric dimension of the Sierpi\'{n}ski product of trees. Notice that for any trees $T_1$ and $T_2$ and any function $f\in H^G$, the Sierpi\'{n}ski product $T_1 \otimes _{f} T_2$ is a tree.

\begin{theorem}
\label{t:thm-tree-1}
If $T_1$ and $T_2$ are trees with $n(T_2)\ge 3$, then
\[
\DimS(T_1,T_2) = n(T_1) \dim(T_2).
\]
\end{theorem}
\proof Let $w$ be a branch vertex of $T_2$, and let $f_w \colon V(T_1) \rightarrow V(T_2)$ be the constant function defined by $f_w(v) = w$ for every vertex $v \in V(T_1)$. The exterior branch vertices in $T_1 \otimes _{f_w} T_2$ are precisely the exterior branch vertices in each of the copies of $T_2$, and so $\ex(T_1 \otimes _{f_w} T_2) = n(T_1)\ex(T_2)$. Moreover, the leaves in $T_1 \otimes _{f_w} T_2$ are precisely the leaves in each of the copies of $T_2$, and so $n_1(T_1 \otimes _{f_w} T_2) = n(T_1)n_1(T_2)$. Therefore, as $T_1 \otimes _{f_w} T_2$ is a tree, by~\eqref{bounDimUsandoExterior} we have
\[
\begin{array}{lcl}
\dim(T_1 \otimes _{f_w} T_2) & = & n_1(T_1 \otimes _{f_w} T_2) - \ex(T_1 \otimes _{f_w} T_2) \1 \\
& = & n(T_1)n_1(T_2) - n(T_1)\ex(T_2) \1 \\
& = & n(T_1)(n_1(T_2) - \ex(T_2)) \1 \\
& = &  n(T_1) \dim(T_2),
\end{array}
\]
implying that
\begin{equation}
\label{lower-bd-1}
\DimS(T_1,T_2) \ge \dim(T_1 \otimes _{f_w} T_2) = n(T_1) \dim(T_2).
\end{equation}

We next show that $\DimS(T_1,T_2) \le n(T_1) \dim(T_2)$. Let $f \colon V(T_1) \rightarrow V(T_2)$ be an arbitrary function. It suffices for us to show that
\[
\dim(T_1 \otimes _{f} T_2) \le n(T_1) \dim(T_2).
\]
Let us first consider that $T_2$ is a path $P_n$ with $n\ge 3$. Hence, we here indeed need to prove that $\dim(T_1 \otimes _{f} P_n) \le n(T_1)$ since $\dim(P_n)=1$. Suppose to the contrary that $\dim(T_1 \otimes _{f} P_n)>n(T_1)$. Thus, from \eqref{bounDimUsandoExterior} we have that
$n_1(T_1 \otimes _{f} P_n)-\ex(T_1 \otimes _{f} P_n)=\dim(T_1 \otimes _{f} P_n)>n(T_1)$. Hence,
\[
\ex(T_1 \otimes _{f} P_n)<n_1(T_1 \otimes _{f} P_n)-n(T_1)\le 2n(T_1)-n(T_1)=n(T_1),
\]
which means there is a positive integer $k$ such that $\ex(T_1 \otimes _{f} P_n)=n(T_1)-k$. Now, notice that if $T_1 \otimes _{f} P_n$ has $n(T_1)-k$ exterior branch vertices, then there must be at least $k$ copies of $P_n$ in $T_1 \otimes _{f} P_n$ not containing any exterior branch vertex. This situation can only happen when the connecting edges of $T_1 \otimes _{f} P_n$ in such copies of $P_n$ are incident with at least one leaf of each of these copies. Consequently, we deduce that $n_1(T_1 \otimes _{f} P_n)\le 2n(T_1)-k$. Therefore, by \eqref{bounDimUsandoExterior} we have
\[
\begin{array}{lcl}
\dim(T_1 \otimes _{f} P_n) & = & n_1(T_1 \otimes _{f} P_n)-\ex(T_1 \otimes _{f} P_n) \1 \\
& \le & 2n(T_1)-k-(n(T_1)-k) \1 \\
& = & n(T_1),
\end{array}
\]
which is a contradiction with our assumption, and so $\dim(T_1 \otimes _{f} P_n) \le n(T_1)$ as required.

We next consider the case when $T_2$ is not a path. Let $E_c = \{e_1, e_2, \ldots, e_{m(T_1)}\}$ be the set of connecting edges in $T_1 \otimes _{f} T_2$.
We order these connecting edges and define $e_i$ as the $i$th \emph{connecting edge} of $T_1 \otimes _{f} T_2$ for $i \in [m(T_1)]$. We next define forests $X_0, X_1, \ldots, X_{m(T_1)}$ as follows. Let $X_0$ be obtained from the tree $T_1 \otimes _{f} T_2$ by removing the connecting edges in $E_c$. We note that $X_0$ is  the disjoint union of $n(T_1)$ copies of the tree $T_2$. Applying~\eqref{bounDimUsandoExterior} to each component of the forest $X_0$ we have
\[
\begin{array}{lcl}
\dim(X_0) & = & n_1(X_0) - \ex(X_0) \1 \\
& = & n(T_1)n_1(T_2) - n(T_1)\ex(T_2) \1 \\
& = & n(T_1)(n_1(T_2) - \ex(T_2)) \1 \\
& = & n(T_1)\dim(T_2).
\end{array}
\]

We now define the forests $X_1, \ldots, X_{m(T_1)}$ as follows. For $i \in [m(T_1)]$, let $X_i$ be the forest obtained from $X_{i-1}$ by adding the $i$th connecting edge, that is, $X_i = X_{i-1} \cup \{e_i\}$. We note that
\[
T_1 \otimes _{f} T_2 = X_{m(T_1)}.
\]
Applying~\eqref{bounDimUsandoExterior} to each component of the forest $X_i$ we have
\[
\dim(X_i) = n_1(X_i) - \ex(X_i)
\]
for all $i \in \{0,1,\ldots,m(T_1)\}$. Let
\[
\Phi_1(X_i) = n_1(X_{i-1}) - n_1(X_i)\,,
\]
which represents the number of vertices of degree~$1$ in $X_{i-1}$ which are of degree at least~$2$ in $X_i$. Roughly speaking, $\Phi_1(X_i)$ is the number of degree~$1$ vertices in $X_{i-1}$ ``destroyed" by adding the $i$th connecting edges $e_i$ to $X_{i-1}$ when constructing $X_i$. Set further
\[
\Phi_2(X_i) = -[\ex(X_{i-1}) - \ex(X_i)]\,, \2
\]
where $\ex(X_{i-1}) - \ex(X_i)$ is the difference between the number of exterior branch vertices in $X_{i-1}$ and the number of exterior branch vertices in $X_i$. We note that
\[
\begin{array}{lcl}
\dim(X_i) & = & n_1(X_i) - \ex(X_i) \1 \\
& = & (n_1(X_{i-1}) - \Phi_1(X_i)) - (\Phi_2(X_i) + \ex(X_{i-1})) \1 \\
& = & \dim(X_{i-1}) - (\Phi_1(X_i) + \Phi_2(X_i))
\end{array}
\]
for all $i \in [m(T_1)]$. We would like to show that
\begin{equation}
\label{relate}
\Phi_1(X_i) + \Phi_2(X_i) \ge 0
\end{equation}
for all $i \in [m(T_1)]$, which would imply that $\dim(X_i) \le \dim(X_{i-1}) \le \dim(X_0) = n(T_1)\dim(T_2)$, from which we deduce that
\[
\dim(T_1 \otimes _{f} T_2) = \dim(X_{m(T_1)}) \le n(T_1)\dim(T_2).
\]
Hence to prove the theorem, it remains to show that~\eqref{relate} holds. For this purpose, let the $i$th connecting edge $e_i$ join vertices $x_i$ and $y_i$ in $X_{i-1}$ when constructing $X_i$.

Suppose that $x_i$ is neither an internal terminal vertex nor a leaf in $X_{i-1}$. In this case, the vertex $x_i$ contributes~$0$ to both terms $\Phi_1(X_i)$ and $\Phi_2(X_i)$.

Suppose that $x_i$ is an internal terminal vertex in $X_{i-1}$. Let $w_i$ be the exterior branch vertex associated with the vertex~$x_i$ in $X_{i-1}$. In this case, when $e_i$ is added to $X_{i-1}$, the vertex $x_i$ becomes an exterior branch vertex in $X_i$, while the vertex~$w_i$ may no longer be an exterior branch vertex, implying that the vertex $x_i$ contributes~$0$ to the term $\Phi_1(X_i)$, and contributes at least~$1 - 1 = 0$ to the term $\Phi_2(X_i)$.

Suppose that $x_i$ is a leaf in $X_{i-1}$. As before, let $w_i$ be the exterior branch vertex associated with the vertex~$x_i$ in $X_{i-1}$. In this case, when $e_i$ is added to $X_{i-1}$, the leaf $x_i$ in $X_{i-1}$ is not a leaf in $X_i$, and therefore the vertex $x_i$ contributes~$1$ to the term $\Phi_1(X_i)$. Moreover, the effect of adding $e_i$ is that the vertex~$w_i$ may no longer be an exterior branch vertex, implying that the vertex $x_i$ contributes at least~$-1$ to the term $\Phi_2(X_i)$. In all of the above three cases, the contribution of $x_i$ to $\Phi_1(X_i) + \Phi_2(X_i)$ is at least~$0$. Analogous arguments hold for the vertex $y_i$, showing that the contribution of $y_i$ to $\Phi_1(X_i) + \Phi_2(X_i)$ is at least~$0$. Therefore,~\eqref{relate} holds. This completes the proof of Theorem~\ref{t:thm-tree-1}.~\QED

\subsection{Sierpi\'{n}ski metric dimension in trees}

In this section we study the Sierpi\'{n}ski metric dimension of the Sierpi\'{n}ski product of trees. The Sierpi\'{n}ski metric dimension of two paths is given by the following result.

\begin{proposition}
\label{t:thm-mdim-paths}
If $T_1$ and $T_2$ are both paths, then $\dimS(T_1,T_2) = 1$.
\end{proposition}
\proof Let $T_1 = P_n$ and let $T_2 = P_m$. If $n = 1$ or $m = 1$, then $T_1 \otimes _{f} T_2$ is a path, and so by Corollary~\ref{known-mdim-tree}, $\dimS(T_1,T_2) = 1$. Hence we may assume that $n \ge 2$ and $m \ge 2$, for otherwise the result is immediate. Let the path $T_2$ be an $x,y$-path that starts at vertex~$x$ and ends at vertex~$y$, and let $T_1$ be the path $v_1 v_2 \ldots v_n$. Let $f \colon V(T_1) \rightarrow V(T_2)$ be the function defined by
\[
f(v_i) =
\left\{
\begin{array}{ll}
x; & \mbox{$i \, (\modo \, 4) \in \{1,2\}$}, \1 \\
y; & \mbox{otherwise}.
\end{array}
\right.
\]
for all $i \in [n]$. In this case the Sierpi\'{n}ski product $T_1 \otimes _{f} T_2$ is a path $P_{nm}$, and so by Corollary~\ref{known-mdim-tree}, $\dimS(T_1,T_2) \le \dim(T_1 \otimes _{f} T_2) = 1$. Consequently, $\dimS(T_1,T_2) = 1$.~\QED

\medskip
In view of Proposition~\ref{t:thm-mdim-paths} it is only of interest to study the Sierpi\'{n}ski product of two trees with at least one of the trees not a path. In this case, we shall establish the following lower bound on the Sierpi\'{n}ski metric dimension, where we use the notation $d_{T}(v)$ to represent the \emph{degree} of a vertex $v$ in $T$.

\begin{lemma}
\label{t:thm-tree-2}
If $T_1$ and $T_2$ are trees, where $T_2$ is not a path, then
\[
\dimS(T_1,T_2) \ge \sum_{v \in V(T_1)} \max \{ 0,\dim(T_2) - d_{T_1}(v) \}.
\]
\end{lemma}
\proof Let $f \colon V(T_1) \rightarrow V(T_2)$ be an arbitrary function. Recall that for each vertex $v \in V(T_1)$, $vT_2$ denotes the copy of the tree $T_2$ in $T_1 \otimes _f T_2$ corresponding to the vertex $v$. We let $C_v$ be the set of vertices in $vT_2$ that are connecting vertices in $T_1 \otimes _{f} T_2$. Thus, each vertex in $C_v$ is incident with a connecting edge in $T_1 \otimes _f T_2$ that joins that vertex to a vertex in a copy of $T_2$ different from $vT_2$. We note that
\[
|C_v| \le d_{T_1}(v) \1
\]
since every edge incident with $v$ in the tree $T_1$ is associated with a connecting edge in $T_1 \otimes _f T_2$ that is incident with a vertex in $vT_2$. Let $B$ be a standard metric basis of $T_1 \otimes _f T_2$, and so $\dim(T_1 \otimes _{f} T_2) = |B|$ and the basis $B$ contains all but one leaf associated with the terminal vertices of degree at least~$2$ in $T_1 \otimes _{f} T_2$. Let $B_v$ be the restriction of $B$ to $vT_2$, that is,
\[
B_v = B \cap V(vT_2)
\]
for every vertex $v \in V(T_1)$. The set $B_v \cup C_v$ is a resolving set in the tree $vT_2$, and so
\[
\dim(T_2) = \dim(vT_2) \le |B_v| + |C_v| \le |B_v| + d_{T_1}(v),
\]
and so $|B_v| \ge \dim(T_2) - d_{T_1}(v)$. As clearly $|B_v| \ge 0$, we get
\[
|B_v| \ge \max\{ 0, \dim(T_2) - d_{T_1}(v) \}
\]
for every vertex $v \in V(T_1)$. Therefore,
\[
\dimS(T_1,T_2) = |B| = \sum_{v \in V(T_1)} |B_v| \ge \sum_{v \in V(T_1)} \max\{ 0, \dim(T_2) - d_{T_1}(v) \}. \1
\]
This establishes the desired lower bound in the statement of the theorem.~\QED

\medskip
Using Lemma~\ref{t:thm-tree-2}, we have the following result.

\begin{theorem}
\label{t:cor-tree-2}
For $n \ge 2$, if $T_1 = P_n$ and $T_2$ is a tree that is not a path, then
\[
\dimS(T_1,T_2) = n(\dim(T_2) - 2) + 2.
\]
\end{theorem}
\proof Let $T_1 = P_n$ and let $T_2$ be a tree that is not a path. By Corollary~\ref{known-mdim-tree}, $\dim(T_2) \ge 2$. By  Lemma~\ref{t:thm-tree-2},
\[
\dimS(T_1,T_2) \ge n(\dim(T_2) - 2) + 2
\]
noting that $T_1$ contains two vertices of degree~$1$ and $n-2$ vertices of degree~$2$. Hence, it suffices for us to show that
\[
\dimS(T_1,T_2) \le n(\dim(T_2) - 2) + 2.
\]

By assumption, $T_2$ is not a path. Hence, $T_2$ contains at least one exterior branch vertex with terminal degree at least~$2$. If $T_2$ contains two distinct exterior branch vertices both with terminal degree at least~$2$, then let $u_1$ and $u_2$ be two selected (terminal) leaves associated with these two exterior branch vertices. If $T_2$ contains only one exterior branch vertex, then $T_2$ is a star or a subdivided star with terminal degree at least~$3$. In this case, let $u_1$ and $u_2$ be two arbitrary leaves in $T_2$. Let $T_1$ be the path $v_1 v_2 \ldots v_n$, and define the function $f \colon V(T_1) \rightarrow V(T_2)$ by
\[
f(v_i) =
\left\{
\begin{array}{ll}
u_1; & \mbox{$i \, (\modo \, 4) \in \{1,2\}$}, \1 \\
u_2; & \mbox{otherwise}.
\end{array}
\right.
\]
for all $i \in [n]$. For notational simplicity, instead of $v_iT_2^i$ we simply write $iT_2$ for the copy of $T_2$ in the Sierpi\'{n}ski product $T_1 \otimes _{f} T_2$ that corresponds to the vertex~$v_i$ for all $i \in [n]$. We note that the copy $1T_2$ and $nT_2$ both contain one less leaf in the product $T_1 \otimes _{f} T_2$, while every copy $iT_2$ where $i \in [n-1] \setminus \{1\}$  contains two fewer leaves in the product $T_1 \otimes _{f} T_2$, implying that
\[
n_1(T_1 \otimes _{f} T_2) = n(T_1)(n_1(T_2) - 2) + 2 = n (n_1(T_2) - 2) + 2.
\]

On the other hand, by our choice of the vertices $u_1$ and $u_2$, the number of exterior branch vertices in each copy $iT_2$ where $i \in [n]$ remains unchanged in the product $T_1 \otimes _{f} T_2$, that is,
\[
\ex(T_1 \otimes _{f} T_2) = n(T_1) \ex(T_2) = n \times \ex(T_2).
\]

Therefore,
\[
\begin{array}{lcl}
\dim(T_1 \otimes _{f} T_2) & = & n_1(T_1 \otimes _{f} T_2) - \ex(T_1 \otimes _{f} T_2) \1 \\
& = & (n(n_1(T_2) - 2) + 2) - n \times \ex(T_2) \1 \\
& = & n(n_1(T_2) - \ex(T_2) - 2) + 2 \1 \\
& = & n (\dim(T_2) - 2) + 2,
\end{array}
\]
completing the proof of Theorem~\ref{t:cor-tree-2}.~\QED

\section{Sierpi\'nski products of cycles}
\label{S:cycles}

In this section we study the Sierpi\'{n}ski metric dimension of the Sierpi\'{n}ski product of cycles, and proved closed formulas for the cases in which at least one of the factors is a triangle.

\begin{theorem}
\label{t:thm-C3-upper}
If $n\ge 3$, then $\DimS(C_n,C_3) = n$.
\end{theorem}
\proof
Let $H = C_3$ and $G = C_n$. Let $V(H) = [3]$, and let the vertices of the cycle $G$ be $g_1, g_2, \ldots, g_n$ in the natural order of adjacencies. Let $f:V(G)\rightarrow V(H)$ be an arbitrary function and consider $G \otimes_{f} H$.

For notational simplicity, let $iH$ denote $g_iH$ for all $i \in [n]$, that is, $iH$ is the $i$th copy of $H$ corresponding to the vertex $g_i$ of $G$. Let $x_iy_{i+1}$ be the connecting edge between $iH$ and $(i+1)H$, $i \in [n]$, where the sum is made mod~$n$. Thus, $x_i = (g_i,f(g_{i+1}))$ and $y_{i+1} = (g_{i+1},f(g_{i}))$. Set further $w_i$ be a vertex in $iH$ different from $x_i$ and $y_i$. Note that if $x_i\ne y_i$, then $V(iH) = \{x_i, y_i, w_i\}$. In case $x_i = y_i$, then let $w_i'$ be the third vertex of $V(iH)$, that is, in this case we have $V(iH) = \{x_i, w_i, w_i'\}$.

Let $Z = \{z_1,z_2, \ldots, z_n\}\subseteq V(G \otimes_{f} H)$ be defined as follows. Consider an arbitrary $iH$, $i\in [n]$. If $x_i \ne y_i$, then let $z_i = y_i$, and if $x_i = y_i$, then let $z_i = w_i$.  Since the vertices $z_i$ are pairwise different, $|Z| = n$. We claim that $Z$ is a resolving set of $G \otimes_{f} H$.

Consider first two vertices from a given $iH$. If $z_i = w_i$, then the vertices $x_i$ and $w_i'$ are distinguished by all the vertices of $Z\setminus \{z_i\}$. If $z_i = y_i$, then $d(w_i,z_{i+1}) = d(x_i,z_{i+1}) + 1$ and therefore $x_i$ and $w_i$ are distinguished. Consider now two vertices $u\in V(iH)$ and $v\in V(jH)$, where $i\ne j$. Assume first that $z_i = w_i$. Then $d(u,z_i) \le 1$ while $d(v,z_i) > 1$ and so $u$ and $v$ are distinguished by $z_i$. The case when $z_j = w_j$ is treated exactly the same. Hence it remains to consider the case when $z_i = y_i$ and $z_j = y_j$. If $|i-j| > 1$, then as before $d(u,z_i) < d(v,z_i)$. Assume finally that $j = i+1$. But then $d(v,z_i) > d(u,z_i)$ and so again $u$ and $v$ are distinguished.

We have thus proved that $Z$ is a resolving set and hence $\DimS(C_n,C_3) \le |Z| = n$. To prove that $\DimS(C_n,C_3) \ge n$, consider the constant function $f_1(g_i) = 1$ for every $i\in [n]$. Then in every $iH$ the vertices corresponding to $2$ and $3$ are twins and hence at least one of them must be included in every resolving set. Then $\dim(G \otimes _{f_1} H) \ge n$ and therefore $\DimS(C_n,C_3) \ge n$.~\QED

\medskip
We remark that, in the proof above, the situation in which we consider the constant function $f_1(g_i) = 1$ for every $i\in [n]$, leads to a graph which indeed represents a rooted product graph, and so, the value of its metric dimension can be also deduced from results appearing in \cite{feng-2013,kuziak-2017}.

For a given integer $k\ge 3$, let $F_k$ be the graph obtained as follows. We begin with a cycle $v_0v_1\cdots v_{2k-1}v_0$ and do all computations modulo $2k$. Next we add $k$ isolated vertices $u_0,u_2,u_4,\dots,u_{2(k-1)}$ and the edges $u_{2i}v_{2i-1},u_{2i}v_{2i}$ for every $0\le i\le k-1$.

\begin{lemma}
\label{lem:all-distances-in-Fk}
If $k\ge 3$, then $\dim(F_k)=2$.
\end{lemma}
\proof
We claim that the set $S=\{u_0,u_{2(\lceil k/2\rceil-1)}\}$ is a metric basis for $F_k$. Let $x,y$ be any two distinct vertices of $F_k$ with $x,y\notin S$ (if $x\in S$ or $y\in S$, then they are clearly identified by $S$). If $d(x,u_0)=d(y,u_0)$, then we have either one of the following situations.

\medskip\noindent
{\bf Case 1}: $x,y\in A_1=\{v_0,v_1,\dots, v_{k-1}\}\cup \{u_0,u_2,\dots,u_{2\lceil k/2\rceil}\}$. \\
In such situation, it must happen that (w.l.g.) $x=v_{2i}$ and $y=u_{2i}$ for some $i\in \{1,2,\dots, \lceil k/2\rceil\}$. If $i=\lceil k/2\rceil-1$, then $y\in S$ and the conclusion is clear. In all the other cases we  notice that $v_{2i}$ belongs to the $u_{2i},u_{2(\lceil k/2\rceil-1)}$-geodesic, which means $u_{2(\lceil k/2\rceil-1)}$ identifies $x=v_{2i}$ and $y=u_{2i}$.

\medskip\noindent
{\bf Case 2}: $x,y\in A_2=V(F_k)\setminus A_1$. \\
A similar conclusion as in Case 1 can be deduced, but taking into account that $x=v_{2i-1}$ and $y=u_{2i}$ for some $i\in \{\lceil k/2\rceil + 1,\dots,k-1\}$.

\medskip\noindent
{\bf Case 3}: $x\in A_1$ and $y\in A_2$. \\
Hence, $x=u_i$ or $x=v_i$, and $y=v_{2k-i-1}$ or $y=u_{2k-i}$, for some $i\in \{0,1,\dots,k-1\}$, where if $x=u_i$, then $i$ is even. We consider now the distances between $x,y$ and $u_{2(\lceil k/2\rceil-1)}$. That is:
\[
\begin{array}{lclcl}
  d(v_i,u_{2(\lceil k/2\rceil-1)}) & = & 2(\lceil k/2\rceil-1)-i & = & 2\lceil k/2 \rceil - i - 2, \\
  d(u_i,u_{2(\lceil k/2\rceil-1)}) & = & 2(\lceil k/2\rceil-1)-i+1 & = &  2\lceil k/2\rceil-i-1, \\
  d(v_{2k-i-1},u_{2(\lceil k/2\rceil-1)}) &=& 2k-i-2(\lceil k/2\rceil-1)& = & 2\lfloor k/2\rfloor -i + 2, \\
  d(u_{2k-i},u_{2(\lceil k/2\rceil-1)}) &=& 2k-i-2(\lceil k/2\rceil-1)+1 & = & 2\lfloor k/2\rfloor -i + 3.
\end{array}
\]
Consequently, we now obtain a contradiction from each situation in which we would suppose that $d(x,u_{2(\lceil k/2\rceil-1)})=d(y,u_{2(\lceil k/2\rceil-1)})$, for any possible assumption taken for $x,y$ as considered before. Therefore, $x,y$ are identified by $u_{2(\lceil k/2\rceil-1)}$, and so, $S$ is a resolving set. Since $F_k$ is not a path, then $S$ it is indeed a metric basis, as claimed.~\QED

\begin{theorem}
\label{t:thm-C3-lower}
If $n\ge 3$, then $\dimS(C_n,C_3) = 2$.
\end{theorem}
\proof
Clearly, $\dimS(C_n,C_3) \ge 2$, hence we only need to prove that $\dimS(C_n,C_3) \le 2$. We use the notation from the first two paragraphs of the proof of Theorem~\ref{t:thm-C3-upper}. In particular, $G = C_n$ and $H=C_3$.

For each $n\ge 3$ we define a function $f_n:V(G) \rightarrow V(H)$ as follows. To do so, let us represent a function $f:V(G) \rightarrow V(H)$ as the vector $(f(g_1), f(g_2), \ldots, f(g_n))$. Then set $f_3 = (1,2,3)$ and $f_4 = (1,1,2,2)$. Let then $n\ge 5$.

\medskip\noindent
{\bf Case 1}: $n$ is odd. \\
Let $B$ be the sequence $3,1,2,3$. Let $f_n$ be defined as follows.
\begin{itemize}
\item If $n = 4(k+1) + 1$, $k\ge 0$, then $f_n = (1,2,3,B,\ldots,B,3,1)$, where $B$ appears $k$ times.
\item If $n = 4(k+1) + 3$, $k\ge 0$, then $f_n = (1,2,3,B,\ldots,B)$, where $B$ appears $k+1$ times.
\end{itemize}

\medskip\noindent
{\bf Case 2}: $n$ is even. \\
Let $C$ be the sequence $2,2,3,3$. Let $f_n$ be defined as follows.
\begin{itemize}
\item If $n = 4(k+1)$, $k\ge 0$, then $f_n = (1,1,C,\ldots,C,2,2)$, where $C$ appears $k$ times.
\item If $n = 4(k+1) + 2$, $k\ge 0$, then $f_n = (1,1,C,\ldots,C)$, where $C$ appears $k+1$ times.
\end{itemize}

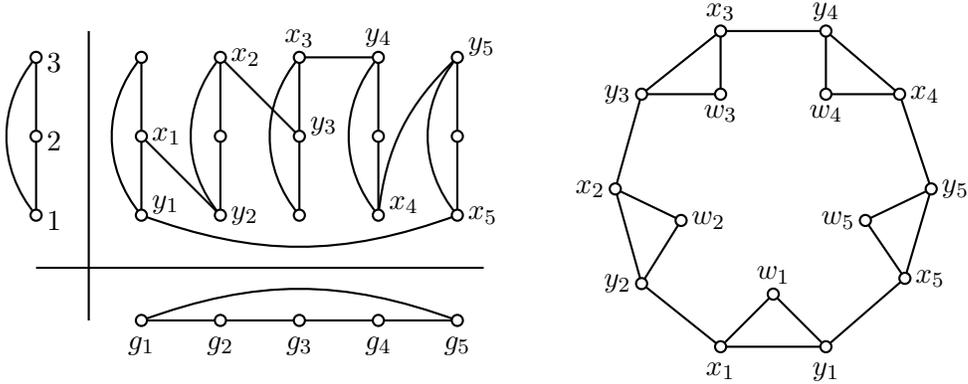
\begin{figure}[ht!]
\begin{center}
\begin{tikzpicture}[scale=0.7,style=thick]
\tikzstyle{every node}=[draw=none,fill=none]
\def\vr{3pt} 

\begin{scope}[yshift = 0cm, xshift = 0cm]
\path (1,-1) coordinate (g1);
\path (2.5,-1) coordinate (g2);
\path (4,-1) coordinate (g3);
\path (5.5,-1) coordinate (g4);
\path (7,-1) coordinate (g5);

\path (-1,1) coordinate (x1);
\path (-1,2.5) coordinate (x2);
\path (-1,4) coordinate (x3);

\path (1,1) coordinate (g1-1);
\path (1,1) coordinate (g1-1);
\path (1,1) coordinate (g1-1);

\path (1,1) coordinate (g1-1);
\path (1,2.5) coordinate (g1-2);
\path (1,4) coordinate (g1-3);
\path (2.5,1) coordinate (g2-1);
\path (2.5,2.5) coordinate (g2-2);
\path (2.5,4) coordinate (g2-3);
\path (4,1) coordinate (g3-1);
\path (4,2.5) coordinate (g3-2);
\path (4,4) coordinate (g3-3);
\path (5.5,1) coordinate (g4-1);
\path (5.5,2.5) coordinate (g4-2);
\path (5.5,4) coordinate (g4-3);
\path (7,1) coordinate (g5-1);
\path (7,2.5) coordinate (g5-2);
\path (7,4) coordinate (g5-3);
\draw (g1-1) -- (g1-2) -- (g1-3);
\draw (g1-1) to[bend left=40] node {} (g1-3);
\draw (g2-1) -- (g2-2) -- (g2-3);
\draw (g2-1) to[bend left=40] node {} (g2-3);
\draw (g3-1) -- (g3-2) -- (g3-3);
\draw (g3-1) to[bend left=40] node {} (g3-3);
\draw (g4-1) -- (g4-2) -- (g4-3);
\draw (g4-1) to[bend left=40] node {} (g4-3);
\draw (g5-1) -- (g5-2) -- (g5-3);
\draw (g5-1) to[bend left=40] node {} (g5-3);
\draw (g1) to[bend left=20] node {} (g5);
\draw (g1) -- (g5);
\draw (x1) -- (x3);
\draw (x1) to[bend left=40] node {} (x3);
\draw (-1,0) -- (7.5,0);
\draw (0,-1) -- (0,4.5);
\draw (g1-2) -- (g2-1);
\draw (g2-3) -- (g3-2);
\draw (g3-3) -- (g4-3);
\draw (g4-1) to[bend left=20] node {} (g5-3);
\draw (g5-1) to[bend left=20] node {}(g1-1);
\draw (g1-1)  [fill=white] circle (\vr);
\draw (g1-2)  [fill=white] circle (\vr);
\draw (g1-3)  [fill=white] circle (\vr);
\draw (g2-1)  [fill=white] circle (\vr);
\draw (g2-2)  [fill=white] circle (\vr);
\draw (g2-3)  [fill=white] circle (\vr);
\draw (g3-1)  [fill=white] circle (\vr);
\draw (g3-2)  [fill=white] circle (\vr);
\draw (g3-3)  [fill=white] circle (\vr);
\draw (g4-1)  [fill=white] circle (\vr);
\draw (g4-2)  [fill=white] circle (\vr);
\draw (g4-3)  [fill=white] circle (\vr);
\draw (g5-1)  [fill=white] circle (\vr);
\draw (g5-2)  [fill=white] circle (\vr);
\draw (g5-3)  [fill=white] circle (\vr);
\draw (g1)  [fill=white] circle (\vr);
\draw (g2)  [fill=white] circle (\vr);
\draw (g3)  [fill=white] circle (\vr);
\draw (g4)  [fill=white] circle (\vr);
\draw (g5)  [fill=white] circle (\vr);
\draw (x1)  [fill=white] circle (\vr);
\draw (x2)  [fill=white] circle (\vr);
\draw (x3)  [fill=white] circle (\vr);

\draw[below] (g1)++(0.0,-0.1) node {$g_1$};
\draw[below] (g2)++(0.0,-0.1) node {$g_2$};
\draw[below] (g3)++(0.0,-0.1) node {$g_3$};
\draw[below] (g4)++(0.0,-0.1) node {$g_4$};
\draw[below] (g5)++(0.0,-0.1) node {$g_5$};
\draw[right] (x1)++(0.0,-0.1) node {$1$};
\draw[right] (x2)++(0.0,-0.1) node {$2$};
\draw[right] (x3)++(0.0,-0.1) node {$3$};
\draw[right] (g1-2)++(0.0,0.0) node {$x_1$};
\draw[right] (g1-1)++(0.0,0.2) node {$y_1$};
\draw[right] (g2-3)++(0.0,0.0) node {$x_2$};
\draw[right] (g2-1)++(0.0,0.0) node {$y_2$};
\draw[above] (g3-3)++(0.0,0.0) node {$x_3$};
\draw[right] (g3-2)++(0.0,0.2) node {$y_3$};
\draw[above] (g4-3)++(0.0,0.0) node {$y_4$};
\draw[right] (g4-1)++(0.0,0.2) node {$x_4$};
\draw[right] (g5-1)++(0.0,0.0) node {$x_5$};
\draw[right] (g5-3)++(0.0,0.2) node {$y_5$};

\end{scope}

\begin{scope}[yshift = -1.5cm, xshift = 12cm]
\path (2,0) coordinate (y1);
\path (0,0) coordinate (x1);
\path (-1.5,1.2) coordinate (y2);
\path (-2,3) coordinate (x2);
\path (-1.5,4.8) coordinate (y3);
\path (0,6) coordinate (x3);
\path (2,6) coordinate (y4);
\path (3.4,4.8) coordinate (x4);
\path (4,3) coordinate (y5);
\path (3.5,1.3) coordinate (x5);
\path (1,1) coordinate (w1);
\path (-0.75,2.4) coordinate (w2);
\path (0,4.8) coordinate (w3);
\path (2,4.8) coordinate (w4);
\path (2.75,2.4) coordinate (w5);
\draw (y1) -- (x1) -- (y2) -- (x2) -- (y3) -- (x3) -- (y4) -- (x4) -- (y5) -- (x5) -- (y1);
\draw (y1) -- (w1) -- (x1);
\draw (y2) -- (w2) -- (x2);
\draw (y3) -- (w3) -- (x3);
\draw (y4) -- (w4) -- (x4);
\draw (y5) -- (w5) -- (x5);
\draw (y1)  [fill=white] circle (\vr);
\draw (x1)  [fill=white] circle (\vr);
\draw (y2)  [fill=white] circle (\vr);
\draw (x2)  [fill=white] circle (\vr);
\draw (y3)  [fill=white] circle (\vr);
\draw (x3)  [fill=white] circle (\vr);
\draw (y4)  [fill=white] circle (\vr);
\draw (x4)  [fill=white] circle (\vr);
\draw (y5)  [fill=white] circle (\vr);
\draw (x5)  [fill=white] circle (\vr);
\draw (w1)  [fill=white] circle (\vr);
\draw (w2)  [fill=white] circle (\vr);
\draw (w3)  [fill=white] circle (\vr);
\draw (w4)  [fill=white] circle (\vr);
\draw (w5)  [fill=white] circle (\vr);
\draw[below] (y1)++(0.0,-0.1) node {$y_1$};
\draw[below] (x1)++(0.0,-0.1) node {$x_1$};
\draw[left] (y2)++(0.0,0.0) node {$y_2$};
\draw[left] (x2)++(0.0,0.0) node {$x_2$};
\draw[left] (y3)++(0.0,0.0) node {$y_3$};
\draw[above] (x3)++(0.0,0.0) node {$x_3$};
\draw[above] (y4)++(0.0,0.0) node {$y_4$};
\draw[right] (x4)++(0.0,0.0) node {$x_4$};
\draw[right] (y5)++(0.0,0.0) node {$y_5$};
\draw[right] (x5)++(0.0,0.0) node {$x_5$};
\draw[above] (w1)++(0.0,0.0) node {$w_1$};
\draw[right] (w2)++(0.0,0.0) node {$w_2$};
\draw[below] (w3)++(0.0,0.0) node {$w_3$};
\draw[below] (w4)++(0.0,0.0) node {$w_4$};
\draw[left] (w5)++(0.0,0.0) node {$w_5$};

\end{scope}
\end{tikzpicture}
\end{center}
\caption{$C_5 \otimes_{f_5} C_3$, where $f_5 = (1,2,3,3,1)$}
\label{fig:Case-5}
\end{figure}

It is straightforward to verify that for every $n\ge 3$, the Sierpi\'nski product $G \otimes_{f_n} H$ has the following structure. For each $i$ we have $x_i\ne y_i$, and hence $w_i$ is the third vertex from $V(iH)$, see Fig.~\ref{fig:Case-5} where $G \otimes_{f_5} H$ is drawn in two different ways. It is now clear that $G \otimes_{f_n} H \cong F_n$ and hence Lemma~\ref{lem:all-distances-in-Fk} completes the argument.~\QED

\section{Convexity property of Sierpi\'{n}ski products}
\label{S:convex}

In this section, we establish a distance convex property of the Sierpi\'{n}ski product of two graphs. Recall that a subgraph $H$ of a graph $G$ is {\em convex} if whenever $u,v\in V(H)$ and $P$ is a shortest $u,v$-path in $G$, then $P$ lies completely in $H$. 

\begin{theorem}
\label{thm:convex}
If $G$ and $H$ be connected graphs, $f \colon V(G)\rightarrow V(H)$, and $g\in V(G)$, then $gH$ is a convex subgraph of $G \otimes _f H$.
\end{theorem}
\proof
Throughout the proof, let $X = G \otimes _f H$. Suppose on the contrary that there exists vertices $u, v\in V(gH)$ such that $d_{X}(u,v) < d_{gH}(u,v)$. Note that this does not happen in trees, hence in the rest we may assume that $G$ contains  cycles. Suppose now that $u$, $v$, and $gH$ are selected such that $d_{X}(u,v)$ is as small as possible among all such counterexamples. Let $u = (g,h)$, $v = (g,h')$, and let $P$ be a shortest $u,v$-path in $X$. Set further $g_1 = g$.

\begin{figure}[ht!]
\begin{center}
\begin{tikzpicture}[scale=0.8,style=thick]
\tikzstyle{every node}=[draw=none,fill=none]
\def\vr{3pt} 

\begin{scope}[yshift = 0cm, xshift = 0cm]
\path (2,-1) coordinate (g1);
\path (5,-1) coordinate (g2);
\path (8,-1) coordinate (g3);
\path (14,-1) coordinate (gk);
\path (2,11) coordinate (u);
\path (2,4) coordinate (v);
\path (-1,11) coordinate (h);
\path (-1,4) coordinate (h');
\path (5,9) coordinate (u2);
\path (5,7) coordinate (v2);
\path (8,11) coordinate (u3);
\path (8,8) coordinate (v3);
\path (14,10) coordinate (uk);
\path (14,6) coordinate (vk);
\draw (-2,0) -- (15,0);
\draw (0,-2) -- (0,12.5);
\draw (u) -- (u2);
\draw (u) -- (u2);
\draw (v2) -- (u3);
\draw (v3) -- (10,7);
\draw (uk) -- (12,9);
\draw (vk) -- (v);
\draw[line width = 2pt, dashed] (5,9) .. controls (6,8.5) and (4,7.5)..(5,7);
\draw[line width = 2pt, dashed] (8,11) .. controls (9,10) and (7,9)..(8,8);
\draw[line width = 2pt, dashed] (14,10) .. controls (15,9) and (13,7)..(14,6);

\draw[rounded corners] (0.5, -1.7) rectangle (15, -0.3) {};
\draw[rounded corners] (-1.7, 0.3) rectangle (-0.3, 12.5) {};
\draw[rounded corners] (1, 0.3) rectangle (3, 12.5) {};
\draw[rounded corners] (4, 0.3) rectangle (6, 12.5) {};
\draw[rounded corners] (7, 0.3) rectangle (9, 12.5) {};
\draw[rounded corners] (13, 0.3) rectangle (15, 12.5) {};
\draw (g1)  [fill=white] circle (\vr);
\draw (g2)  [fill=white] circle (\vr);
\draw (g3)  [fill=white] circle (\vr);
\draw (gk)  [fill=white] circle (\vr);
\draw (u)  [fill=white] circle (\vr);
\draw (v)  [fill=white] circle (\vr);
\draw (h)  [fill=white] circle (\vr);
\draw (h')  [fill=white] circle (\vr);
\draw (u2)  [fill=white] circle (\vr);
\draw (v2)  [fill=white] circle (\vr);
\draw (u3)  [fill=white] circle (\vr);
\draw (v3)  [fill=white] circle (\vr);
\draw (uk)  [fill=white] circle (\vr);
\draw (vk)  [fill=white] circle (\vr);
\draw[above] (u)++(0.0,0.0) node {\small{$(g_1,f(g_2))$}};
\draw[below] (u)++(0.0,0.0) node {\small{$u$}};
\draw[below] (h)++(0.0,0.0) node {\small{$h$}};
\draw[below] (h')++(0.0,0.0) node {\small{$h'$}};
\draw[above] (v)++(0.0,0.0) node {\small{$(g_1,f(g_k))$}};
\draw[below] (v)++(0.0,0.0) node {\small{$v$}};
\draw[above] (u2)++(0.0,0.4) node {\small{$(g_2,f(g_1))$}};
\draw[below] (v2)++(0.0,0.0) node {\small{$(g_2,f(g_3))$}};
\draw[above] (u3)++(0.0,0.0) node {\small{$(g_3,f(g_2))$}};
\draw[below] (v3)++(0.0,-0.3) node {\small{$(g_3,f(g_4))$}};
\draw[above] (uk)++(0.0,0.0) node {\footnotesize{$(g_k,f(g_{k-1}))$}};
\draw[below] (vk)++(0.0,0.0) node {\small{$(g_k,f(g_1))$}};
\draw[below] (g1)++(0.0,0.0) node {\small{$g=g_1$}};
\draw[below] (g1)++(0.0,0.0) node {\small{$g=g_1$}};
\draw[below] (g2)++(0.0,0.0) node {\small{$g_2$}};
\draw[below] (g3)++(0.0,0.0) node {\small{$g_3$}};
\draw[below] (gk)++(0.0,0.0) node {\small{$g_k$}};
\draw (0.8,-1) node {$G$};
\draw (-1,1) node {$H$};
\draw (4.5,8) node {$P_2$};
\draw (8.5,9.5) node {$P_3$};
\draw (14.5,8) node {$P_k$};
\end{scope}

\end{tikzpicture}
\end{center}
\caption{The shape of $P$}
\label{fig:shape-of-P}
\end{figure}
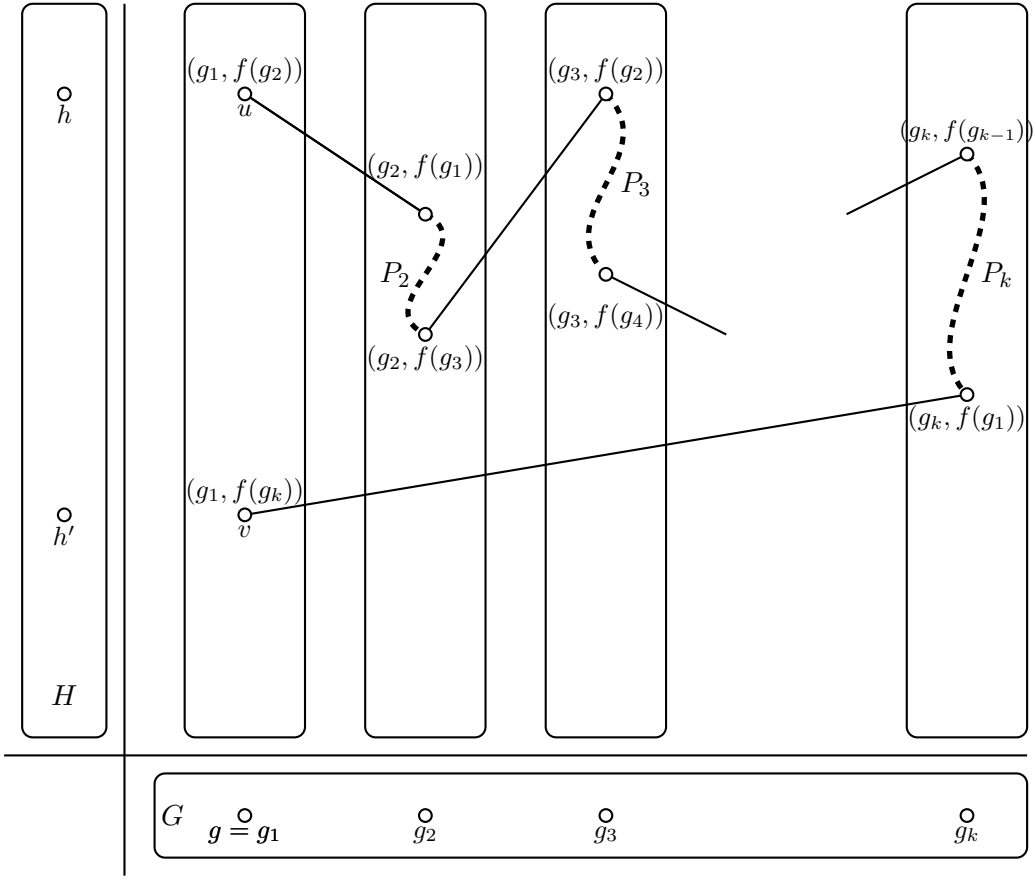

\medskip\noindent
{\bf Claim}. The shape of the path $P$ is as follows. Let $g_1, \ldots, g_k$, $k\ge 2$, be the vertices of $G$ ordered such that $P$ passes through $g_1H, \ldots, g_kH$ in that order. Then $P$ starts with the connecting edge $(g_1,f(g_2))(g_2,f(g_1))$, proceeds with a geodesic $P_2$ in $g_2H$ between $(g_2,f(g_1))$ and $(g_2,f(g_3))$, then continuing with the connecting edge $(g_2,f(g_3))(g_3,f(g_2))$, and so on. Finally $P$ arrives at $g_kH$, proceeds along a geodesic in $G_k$ between $(g_k,f(g_{k-1}))$ and $(g_k,f(g_1))$, and ends with the connecting edge $(g_k,f(g_1))(g_1,f(g_k))$, where $f(g_k) = h'$. See Fig.~\ref{fig:shape-of-P}.

Note first that $k\ge 3$ because $k=2$ would imply that there are two connecting edges between $g_1H$ and $g_2H$.

We next show that the vertices $g_1, \ldots, g_k$ are pairwise different. Suppose on the contrary that there exist $i$ and $j$, such that $2\le i < j < k$ and $g_{j+1} = g_i$. Let $P'$ be the subpath of $P$ between the vertices $u' = (g_i,f(g_{i+1}))$ and $v' = (g_{j+1},f(g_{j})) = (g_{i},f(g_{j}))$. As $P$ is a geodesic in $X$, Bellman's principle of optimality implies that $P'$ is also a geodesic. In addition to its connecting edges, $P'$ contains geodesics $P_i, P_{i+1}, \ldots, P_j$, which are respectively, projected onto $H$, geodesics
\begin{itemize}
\item between $f(g_{i-1})$ and $f(g_{i+1})$,
\item between $f(g_{i})$ and $f(g_{i+2})$,\\
$\vdots$
\item between $f(g_{j-2})$ and $f(g_{j})$, and
\item between $f(g_{j-1})$ and $f(g_{i})$.
\end{itemize}
Suppose first that $j - i$ is odd. Then we have
\begin{align*}
|P'| & > d_H(f(g_{i-1}), f(g_{i+1})) + d_H(f(g_{i+1}), f(g_{i+3})) + \cdots + d_H((g_{j-2}), f(g_{j})) \\
& \ge d_H(f(g_{i-1}), f(g_{j})) \\
& = d_{g_iH}(u',v')\,.
\end{align*}
This is a contradiction with the selection of $u$ and $v$ as a minimal counterexample.
Suppose second that $j - i$ is even. Then we have
\begin{align*}
|P'| & > [d_H(f(g_{i-1}), f(g_{i+1})) + d_H(f(g_{i+1}), f(g_{i+3})) + \cdots + d_H((g_{j-1}), f(g_{i}))] + \\
& \quad\ [d_H(f(g_{i}), f(g_{i+2})) + d_H(f(g_{i+2}), f(g_{i+4})) + \cdots + d_H((g_{j-2}), f(g_{j}))] \\
& \ge d_H(f(g_{i-1}), f(g_{i})) + d_H(f(g_{i}), f(g_{j}))\\
& \ge d_H(f(g_{i-1}), f(g_{j}))\\
& = d_{g_iH}(u',v')\,.
\end{align*}
Hence we get the same contradiction as in the previous case.

We have thus proved that the vertices $g_1, \ldots, g_k$ are pairwise different. To complete the proof of the claim, we need to verify that $P$ starts and ends with a connecting edge. Suppose on the contrary that $P$ starts with a subpath in $g_1H$ from $u$ to $w$ and then proceed along the connecting edge between $g_1H$ and $g_2H$. By the minimality assumption on $u$ and $v$ we have $d_X(u,w) = d_{g_1H}(u,w)$. We now have:
\begin{align*}
d_X(u,w) + d_X(w,v) & = d_X(u,v) \\
& < d_{g_1H}(u,v) \\
& \le d_{g_1H}(u,w) + d_{g_1H}(w,v) \\
& = d_X(u,w) + d_{g_1H}(w,v)
\end{align*}
which yields $d_X(w,v) < d_{g_1H}(w,v)$. This contradiction proves that $P$ indeed starts with a connecting edge. A parallel argument yields that $P$ also ends with a connecting edge. This proves the claim.

To conclude the proof, let $P_2, P_3, \ldots, P_k$ be the sections of $P$ restricted to $g_2H, g_3H, \ldots, g_kH$, respectively. Then we proceed similarly as we did for the subpath $P'$ above. More precisely, if $k-1$ is odd, then
\begin{align*}
|P| & > d_H(f(g_{2}), f(g_{4})) + d_H(f(g_{4}), f(g_{6})) + \cdots + d_H((g_{k-2}), f(g_{k})) \\
& \ge d_H(f(g_{2}), f(g_{k})) \\
& = d_{g_1H}(u,v)\,.
\end{align*}
And if $k-1$ is even, then
\begin{align*}
|P'| & > [d_H(f(g_{1}), f(g_{3})) + d_H(f(g_{3}), f(g_{5})) + \cdots + d_H((g_{k-2}), f(g_{k}))] + \\
& \quad\ [d_H(f(g_{2}), f(g_{4})) + d_H(f(g_{4}), f(g_{6})) + \cdots + d_H((g_{k-1}), f(g_{1}))] \\
& \ge d_H(f(g_{1}), f(g_{k})) + d_H(f(g_{2}), f(g_{1}))\\
& \ge d_H(f(g_{2}), f(g_{k}))\\
& = d_{g_1H}(u,v)\,.
\end{align*}
This final contradiction proves the theorem.~\QED

\section{Concluding remarks}
\label{S:remark}

In Theorem~\ref{t:cor-tree-2} we determined $\dimS(P_n,T)$ where $T$ is an arbitrary tree different from a path. It remains to determine $\dimS(T',T)$ where $T'$ and $T$ are arbitrary trees.

In Theorems~\ref{t:thm-C3-upper} and~\ref{t:thm-C3-lower} we determined  $\DimS(C_n,C_3)$ and $\dimS(C_n,C_3)$ for all $n \ge 3$. It remains to determine $\DimS(C_n,C_m)$ and $\dimS(C_n,C_m)$ for all $n \ge 3$ and $m \ge 4$.

It would be interesting to determine $\DimS(G,H)$ and $\dimS(G,H)$ for other classes of graphs $G$ and $H$.

\section*{Acknowledgments}

Research of Michael Henning was supported in part by the University of Johannesburg, and obtained during his sabbatical visit at the University of Ljubljana. Sandi Klav\v{z}ar was supported by the Slovenian Research Agency (ARRS) under the grants P1-0297, J1-2452, and N1-0285. Ismael G. Yero has been partially supported by the Spanish Ministry of Science and Innovation through the grant PID2019-105824GB-I00. Moreover, this investigation was developed while this author (Ismael G. Yero) was visiting the University of Ljubljana, Slovenia, supported by ``Ministerio de Educaci\'on, Cultura y Deporte'', Spain, under the ``Jos\'e Castillejo'' program for young researchers (reference number: CAS21/00100).

\section*{Declaration of interests}

The authors declare that they have no known competing financial interests or personal relationships that could have appeared to influence the work reported in this paper.

\section*{Data availability}

Our manuscript has no associated data.

\end{document}